\documentclass{amsart}
\usepackage{amscd}

\theoremstyle{plain}
\newtheorem{thm}{Theorem}[section]
\newtheorem{lem}[thm]{Lemma}
\newtheorem{cor}[thm]{Corollary}

\theoremstyle{definition}
\newtheorem{dfn}[thm]{Definition}
\newtheorem{exa}[thm]{Example}

\newcommand{\C}{\mathbb{C}}

\newcommand{\Z}{\mathbb{Z}}
\newcommand{\N}{\mathbb{N}}

\DeclareMathOperator{\GL}{GL}
\DeclareMathOperator{\tr}{tr}
\DeclareMathOperator{\Hom}{Hom}
\DeclareMathOperator{\Ker}{Ker}
\DeclareMathOperator{\Ima}{Im}
\DeclareMathOperator{\Aut}{Aut}
\DeclareMathOperator{\End}{End}
\DeclareMathOperator{\Ext}{Ext}
\DeclareMathOperator{\Rep}{Rep}
\DeclareMathOperator{\Mat}{Mat}
\DeclareMathOperator{\rank}{rank}

\DeclareMathOperator{\dbslash}{/\!\!/}
\DeclareMathOperator*{\rightoverleft}{\parbox{2em}{\centerline{$\longrightarrow$}\vskip
-6pt\centerline{$\longleftarrow$}}}

\begin{document}
\title[Multiplicative preprojective algebras]{Multiplicative
preprojective algebras, middle convolution and the Deligne-Simpson problem}

\author{William Crawley-Boevey}
\address{Department of Pure Mathematics, University of Leeds, Leeds LS2 9JT, UK}
\email{w.crawley-boevey@leeds.ac.uk}

\author{Peter Shaw}
\address{Department of Pure Mathematics, University of Leeds, Leeds LS2 9JT, UK}
\email{pshaw@maths.leeds.ac.uk}

%15  (1940-now) Linear and multilinear algebra; matrix theory
%15A  (1973-1999) Linear and multilinear algebra; matrix theory
%15A24  (1973-now) Matrix equations and identities

%16  (1959-now) Associative rings and algebras [For the commutative case, see 13-XX]
%16G  (1991-now) Representation theory of rings and algebras
%16G20  (1991-now) Representations of quivers and partially ordered sets

%34  (1940-now) Ordinary differential equations
%34M  (2000-now) Differential equations in the complex domain [See also 30Dxx, 32G34]
%34M50  (2000-now) Inverse problems (Riemann-Hilbert, inverse differential Galois, etc.)


\begin{abstract}
We introduce a family of algebras which are multiplicative analogues of
preprojective algebras, and their deformations, as introduced by
M.~P.~Holland and the first author. We show that these algebras provide
a natural setting for the `middle convolution' operation
introduced by N.~M.~Katz in his book `Rigid local systems', and put in an
algebraic setting by M.~Dettweiler and S.~Reiter, and by H.~V\"olklein.
We prove a homological formula relating the dimensions of Hom and
Ext spaces, study varieties of representations of multiplicative
preprojective algebras, and use these results to study simple representations.
We apply this work to the Deligne-Simpson problem, obtaining a sufficient
(and conjecturally necessary) condition for the existence of an irreducible
solution to the equation
$A_1A_2\dots A_k=1$ with the $A_i$ in prescribed conjugacy classes in
$\mathrm{GL}_n(\mathbb{C})$.
\end{abstract}
\maketitle

\section{Introduction}
Given conjugacy classes $C_1,\dots,C_k$ in $\GL_n(\C)$, we consider the
following problem: determine whether or not one can find an irreducible solution
to the equation
\begin{equation}
\label{eq1}
A_1 A_2 \dots A_k = 1
\end{equation}
with $A_i\in C_i$. Here `irreducible' means that the $A_i$ have no
common invariant subspace. This problem has been studied before, in
particular by Deligne and Simpson \cite{Simp}, by Katz
\cite{Katz} and by Kostov \cite{Kostovmon,KostovCR,Kostovaspects,Kostovsurvey},
who calls it the `Deligne-Simpson problem'. In \cite{CBipb} the first author
gave a conjectural answer, and in this paper we prove
one direction of the conjecture. Before describing it we
recall the root system associated to a quiver.

Let $Q$ be a quiver with vertex set $I$. For each arrow $a\in Q$
we denote its head and tail vertices by $h(a),t(a)\in I$. Both $I$
and $Q$ are assumed to be finite. We say that a vertex $v$
is \emph{loopfree} if there is no arrow $a$ with $h(a)=t(a)=v$. We
often call elements of $\Z^I$ \emph{dimension vectors}. There is a
quadratic form $q$ on $\Z^I$ given by
\[
q(\alpha) = \sum_{v\in I} \alpha_v^2 - \sum_{a\in Q} \alpha_{h(a)}\alpha_{t(a)}.
\]
Let $(-,-)$ be the corresponding symmetric bilinear form with $(\alpha,\alpha)=2q(\alpha)$.
For later use we also define $p(\alpha)=1-q(\alpha)$.
For a loopfree vertex $v$, the reflection $s_v : \Z^I \to \Z^I$, is defined by
$s_v(\alpha) = \alpha - (\alpha,\epsilon_v) \epsilon_v$,
the \emph{Weyl group} $W$ is the
subgroup of $\Aut(\Z^I)$ generated by the $s_v$, and the
\emph{real roots} are the images under elements of $W$ of the
coordinate vectors $\epsilon_v$ at loopfree vertices $v$.
The \emph{fundamental region} consists of the
nonzero elements $\alpha\in\N^I$ which have connected support and
$(\alpha,\epsilon_v)\le 0$ for all $v$; its closure under
the action of $W$ and change of sign is, by definition, the set of
\emph{imaginary roots}. Note that $p(\alpha)=0$ for real roots,
and $p(\alpha) > 0$ for imaginary roots. Recall that any root is
\emph{positive} ($\alpha\in\N^I$), or \emph{negative}
($-\alpha\in\N^I$).

To fix conjugacy classes $C_1,\dots,C_k$ in $\GL_n(\C)$ we fix a collection of positive integers $w=(w_1,\dots,w_k)$,
and elements $\xi_{ij}\in \C^*$ ($1\le i\le k$, $1\le j\le w_i$) with
\begin{equation}\label{e:xieq}
(A_i - \xi_{i1} 1) (A_i - \xi_{i2} 1) \dots (A_i - \xi_{i,w_i} 1) = 0
\end{equation}
for $A_i\in C_i$. Clearly, if one wishes, one can take $w_i$ to be
the degree of the minimal polynomial of $A_i$, and
$\xi_{i1},\dots,\xi_{i,w_i}$ to be its roots, counted with
the appropriate multiplicity. The conjugacy class $C_i$ is then determined by the
ranks of the partial products
\[
\alpha_{ij} = \rank (A_i - \xi_{i1}1)(A_i - \xi_{i2}1)\dots(A_i - \xi_{ij}1)
\]
for $A_i\in C_i$ and $1\le j\le w_i-1$. Setting $\alpha_0 = n$, we
obtain a dimension vector $\alpha$ for the following quiver $Q_w$
\setlength{\unitlength}{1.5pt}
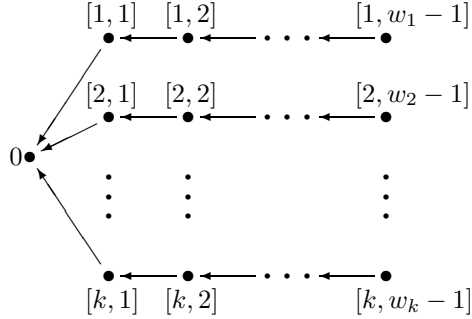
\[
\begin{picture}(110,80)
\put(10,40){\circle*{2.5}} \put(30,10){\circle*{2.5}}
\put(30,50){\circle*{2.5}} \put(30,70){\circle*{2.5}}
\put(50,10){\circle*{2.5}} \put(50,50){\circle*{2.5}}
\put(50,70){\circle*{2.5}} \put(100,10){\circle*{2.5}}
\put(100,50){\circle*{2.5}} \put(100,70){\circle*{2.5}}
\put(28,67){\vector(-2,-3){16}}
\put(27,48.5){\vector(-2,-1){14}}
\put(28,13){\vector(-2,3){16}}
\put(47,10){\vector(-1,0){14}}
\put(47,50){\vector(-1,0){14}}
\put(47,70){\vector(-1,0){14}}
\put(67,10){\vector(-1,0){14}}
\put(67,50){\vector(-1,0){14}}
\put(67,70){\vector(-1,0){14}}
\put(97,10){\vector(-1,0){14}}
\put(97,50){\vector(-1,0){14}}
\put(97,70){\vector(-1,0){14}}
\put(70,10){\circle*{1}}
\put(75,10){\circle*{1}} \put(80,10){\circle*{1}}
\put(70,50){\circle*{1}} \put(75,50){\circle*{1}}
\put(80,50){\circle*{1}} \put(70,70){\circle*{1}}
\put(75,70){\circle*{1}} \put(80,70){\circle*{1}}
\put(30,25){\circle*{1}} \put(30,30){\circle*{1}}
\put(30,35){\circle*{1}} \put(50,25){\circle*{1}}
\put(50,30){\circle*{1}} \put(50,35){\circle*{1}}
\put(100,25){\circle*{1}} \put(100,30){\circle*{1}}
\put(100,35){\circle*{1}} \put(5,38){0} \put(24,2){$[k,1]$}
\put(24,54){$[2,1]$} \put(24,74){$[1,1]$} \put(44,2){$[k,2]$}
\put(44,54){$[2,2]$} \put(44,74){$[1,2]$} \put(92,2){$[k,w_k-1]$}
\put(92,54){$[2,w_2-1]$} \put(92,74){$[1,w_1-1]$}
\end{picture}
\]
with vertex set
$I = \{ 0 \} \cup \{ [i,j] : 1\le i\le k, 1\le j\le w_i-1 \}$.
(For ease of notation, if $\alpha\in\Z^I$, we write
its components as $\alpha_0$ and $\alpha_{ij}$.)
For any $\beta\in\Z^I$ we define
\[
\xi^{[\beta]}
= \prod_{i=1}^k \xi_{i1}^{\beta_0} \prod_{j=1}^{w_i-1} (\xi_{i,j+1}/\xi_{ij})^{\beta_{ij}}
= \prod_{i=1}^k \prod_{j=1}^{w_i} \xi_{ij}^{\beta_{i,j-1} - \beta_{ij}}
\]
using the convention that $\beta_{i0} = \beta_0$ and $\beta_{i,w_i}=0$ for all $i$.

\begin{thm}
\label{t:DSP}
Let $C_1,\dots,C_k$ be conjugacy classes in $\GL_n(\C)$.
Choose $w$ and $\xi$ as above,
and let $\alpha$ be the corresponding dimension vector.
If $\alpha$ is a positive root for $Q_w$,
$\xi^{[\alpha]}=1$, and $p(\alpha) > p(\beta)+p(\gamma)+\dots$ for
any nontrivial decomposition of $\alpha$ as a sum of positive
roots $\alpha = \beta+\gamma+\dots$ with $\xi^{[\beta]} =
\xi^{[\gamma]} = \dots = 1$, then there is an irreducible solution
to $A_1\dots A_k = 1$ with $A_i\in C_i$.
\end{thm}

The proof of this result depends on two ingredients. On the
one hand, we use the result in \cite{CBipb} on solutions to equation
(\ref{eq1}) with the $A_i$ in the closures $\overline{C_i}$ of the conjugacy classes. On
the other hand, we use the properties of a class of algebras
which we now introduce, called `multiplicative preprojective algebras'.

Let $K$ be a field and let $Q$ be a quiver with vertex set $I$.
Recall that the path algebra $KQ$ has as basis the paths
$a_1\dots a_n$ with $a_i\in Q$ and $t(a_i)=h(a_{i+1})$ for
all $i$, and trivial paths $e_v$ ($v\in I$).
Let $\overline{Q}$ be the \emph{double} of $Q$, obtained by
adjoining a reverse arrow $a^*$ for each arrow $a\in Q$.
We extend the operation $a\mapsto a^*$ to an involution on $\overline{Q}$
by defining $(a^*)^* = a$ for $a\in Q$, and we define
$\epsilon(a)=1$ if $a\in Q$ and $\epsilon(a)=-1$ if $a^*\in Q$.

Fix a total ordering $<$ on the set of arrows in
$\overline{Q}$. Given $q\in (K^*)^I$, we consider algebra homomorphisms
$K\overline{Q}\to R$, where $R$ is a $K$-algebra, with the properties that
\begin{gather}
\label{c:inv}
\text{$1 + a a^*$ is invertible in $R$ for all $a\in\overline{Q}$, and}
\\
\label{c:rel}
\text{$\displaystyle\prod_{a\in\overline{Q}} (1+a a^*)^{\epsilon(a)} = \sum_{v\in I} q_v e_v$ in $R$.}
\end{gather}
Here the product has to be taken in the right order,
so if the arrows in $\overline{Q}$ are $a_1<a_2<\dots <a_n$,
then condition (\ref{c:rel}) is that
\[
(1+a_1 a_1^*)^{\epsilon(a_1)} (1+a_2 a_2^*)^{\epsilon(a_2)}\dots
(1+a_n a_n^*)^{\epsilon(a_n)} = \sum_{v\in I} q_v e_v.
\]
It is easy to see that there is a universal such homomorphism, unique up to isomorphism.
It can be constructed by adjoining inverses for each of the elements in (\ref{c:inv}),
and then factoring out the relation (\ref{c:rel}).

\begin{dfn}
We denote the universal homomorphism satisfying (\ref{c:inv}) and (\ref{c:rel}) by $K\overline{Q}\to\Lambda^q$.
The algebra $\Lambda^q$ is called a \emph{multiplicative preprojective algebra}.
If we need to specify the quiver $Q$ and the ordering on $\overline{Q}$ we denote it $\Lambda^q(Q,<)$.
\end{dfn}

Compare this with the definition of the deformed preprojective algebra \cite{CBH}.

\begin{exa}
If $Q$ consists of a single vertex $v$ and a single loop $a$,
then $q\in (K^*)^I$ can be identified with a single element $q\in K^*$.
If the ordering is $a < a^*$,
then $\Lambda^q$ is the algebra given by generators $a$, $a^*$, $\iota_a$, $\iota_{a^*}$
and relations
\begin{gather*}
\iota_a (1+aa^*) = (1+aa^*) \iota_a = 1, \\
\iota_{a^*} (1+a^*a) = (1+a^*a) \iota_{a^*} = 1, \\
(1+aa^*)\iota_{a^*}=q\, 1.
\end{gather*}
\noindent
If $q=1$ this can be rewritten as a localized polynomial algebra
\[
\Lambda^1 \cong K[x,y,(1+xy)^{-1}],
\]
while if $q\neq 1$ then $\Lambda^q$ is isomorphic to
a localized first quantized Weyl algebra,
\[
\Lambda^q \cong B^q_1,
\]
as discussed in \cite{Jordsimp} and the references therein.
\end{exa}

In Section~\ref{s:indep} we prove the following.

\begin{thm}
\label{t:indeporor}
Up to isomorphism, $\Lambda^q$ doesn't depend on the orientation of $Q$
or the chosen ordering on $\overline{Q}$.
\end{thm}

Recall that the category of representations of $KQ$, that is, left $KQ$-modules,
is equivalent to the category of representations of $Q$ by means of vector spaces $X_v$
for each vertex $v$ and linear maps $X_a:X_{t(a)}\to X_{h(a)}$ for each arrow $a$.
The representations of $\Lambda^q$ can be identified with representations
$X$ of $\overline{Q}$ which satisfy
\begin{gather}
\label{c:invtwoprime}
\text{$1_{X_{h(a)}}+X_a X_{a^*}$ is an invertible endomorphism
of $X_{h(a)}$ for all $a\in \overline{Q}$, and}
\\
\label{c:reltwoprime}
\text{$\prod_{\substack{a\in\overline{Q} \\ h(a)=v}}
(1_{X_{h(a)}}+X_{a} X_{a^*})^{\epsilon(a)} = q_v 1_{X_v}$
for all $v\in I$.}
\end{gather}
(To see this, note that the multiplicative preprojective algebra is also universal for
homomorphisms $K\overline{Q}\to R$ satisfying
\begin{gather*}
%\label{c:invprime}
\text{$e_{h(a)} + a a^*$ is invertible in $e_{h(a)}Re_{h(a)}$ for all $a\in\overline{Q}$, and}
\\
%\label{c:relprime}
\text{$\displaystyle\prod_{\substack{a\in\overline{Q} \\ h(a)=v}} (e_v+a a^*)^{\epsilon(a)} = q_v e_v$ in $e_v R e_v$ for all $v\in I$.}
\end{gather*}
Here we use that if $e$ is an idempotent in an algebra $R$ and $x\in eRe$,
then $x+(1-e)$ is invertible in $R$ if and only if $x$ is invertible in $eRe$.)

The \emph{dimension vector} of a finite-dimensional representation
$X$ is the element $\alpha = \underline{\dim} X$ in $\N^I$,
defined by $\alpha_v = \dim e_v X = \dim X_{v}$. For any $\alpha\in\Z^I$, we define
\[
q^\alpha = \prod_{v\in I} q_v^{\alpha_v}.
\]
By considering the determinants of the relations (\ref{c:reltwoprime}), and using the
fact that
\begin{equation}
\label{e:detformula}
\det(1_V+\theta\phi)=\det(1_U+\phi\theta)
\end{equation}
for linear maps $\theta:U\to V$ and $\phi:V\to U$,
we clearly have the following.

\begin{lem}
\label{l:qalphaone}
If $\Lambda^q$ has a representation of dimension vector $\alpha$, then $q^\alpha=1$.
\end{lem}

In Section~\ref{s:homform} we prove the following result.

\begin{thm}
\label{t:homform}
If $X$ and $Y$ are finite-dimensional representations of $\Lambda^q$ then
\[
\dim\Ext^1_{\Lambda^q}(X,Y) = \dim\Hom_{\Lambda^q}(X,Y)+\dim\Hom_{\Lambda^q}(Y,X)-(\underline{\dim}X,\underline{\dim}Y).
\]
\end{thm}

It follows that $\dim\Ext^1_{\Lambda^q}(X,Y) = \dim\Ext^1_{\Lambda^q}(Y,X)$.

One of the reasons for introducing multiplicative preprojective
algebras is to better understand the middle convolution operation
of Katz~\cite{Katz} and its algebraic versions due to Dettweiler
and Reiter~\cite{DR}, and to V\"olklein~\cite{V}.
In Section~\ref{s:MC} we adapt \cite{DR} as follows.
If $v$ is a loopfree vertex in $Q$, define
\[
u_v : (K^*)^I \to (K^*)^I,
\quad
u_v(q)_w = q_v^{-(\epsilon_v,\epsilon_w)} q_w.
\]
Observe that $u_v$ is a multiplicative dual to the reflection $s_v$
on dimension vectors, in the sense that
\begin{equation}
\label{e:usduality}
(u_v(q))^\alpha = q^{s_v(\alpha)}
\end{equation}
for all $\alpha\in\Z^I$.

\begin{thm}
\label{t:MC}
If $v$ is a loopfree vertex and $q_v\neq 1$, then there is an equivalence $F_q$
from the category of representations of $\Lambda^q$ to the category of representations of $\Lambda^{u_v(q)}$.
It acts on dimension vectors as the reflection $s_v$. The inverse equivalence is $F_{u_v(q)}$.
\end{thm}

Observe that in this setting one has an equivalence on the whole of the category of
representations, in contrast to \cite{DR} where it was necessary to
impose certain conditions, denoted (*) and (**).

In Section~\ref{s:simp} we use middle convolution to prove the following result.
(Of course the last part follows from Schur's Lemma if the base field $K$ is algebraically closed.)

\begin{thm}
\label{t:simps}
If there is a simple representation $X$ of $\Lambda^q$ of dimension $\alpha$,
then $\alpha$ is a positive root for $Q$.
If in addition $\alpha$ is a real root, then $\End_{\Lambda^q}(X) = K$.
\end{thm}

We say that a finite-dimensional simple representation $X$ of $\Lambda^q$
is \emph{rigid} if it has no self-extensions, that is, $\Ext^1_{\Lambda^q}(X,X)=0$.
In view of Theorems~\ref{t:homform} and~\ref{t:simps},
it is equivalent that $\underline{\dim}X$ is a real root.
Note that if $X$ is a rigid simple representation of
dimension $\alpha$, then any other representation $Y$ of
dimension $\alpha$ must be isomorphic to $X$, for
Theorem~\ref{t:homform} guarantees the existence of
a nonzero homomorphism $X\to Y$ or $Y\to X$, but since
$X$ is simple any such homomorphism must be an isomorphism.
(Compare with \cite[Lemma 6]{Simp}.)
Middle convolution gives the following characterization of the
possible dimension vectors of rigid simple representations.

\begin{thm}
\label{t:rigidsimp}
There is a rigid simple representation of $\Lambda^q$ of dimension vector $\alpha$
if and only if $\alpha$ is a positive real root, $q^\alpha=1$, and there is
no decomposition $\alpha=\beta+\gamma+\dots$ as a sum of two or more positive
roots with $q^\beta=q^\gamma=\dots=1$.
\end{thm}

Now assume that the field $K$ is algebraically closed.
Recall that the variety of representations of a quiver $Q$ of dimension vector $\alpha$
is the space
\[
\Rep(Q,\alpha) = \prod_{a\in Q} \Mat(\alpha_{h(a)}\times \alpha_{t(a)},K),
\]
and that isomorphism classes correspond to orbits of the algebraic group
\[
\GL(\alpha) = \prod_{v\in I} \GL(\alpha_v,K).
\]
We define $\Rep(\Lambda^q,\alpha)$ to be the subset of
$\Rep(\overline{Q},\alpha)$ consisting of the representations
which satisfy (\ref{c:invtwoprime}) and (\ref{c:reltwoprime}).
It is clear that this is a locally closed subset of $\Rep(\overline{Q},\alpha)$,
so a variety. In Section~\ref{s:var} we prove the following result.

\begin{thm}
\label{t:var}
$\Rep(\Lambda^q,\alpha)$ is an affine variety, and every irreducible component
has dimension at least $g+2p(\alpha)$, where $g = -1+\sum_{v\in I}\alpha_v^2$.
Moreover, the representations $X$ with trivial endomorphism algebra, $\End(X)=K$,
form an open subset of $\Rep(\Lambda^q,\alpha)$ which if nonempty is smooth
of dimension $g+2p(\alpha)$.
\end{thm}

This theorem shows that if $X$ is a non-rigid finite-dimensional simple representation of $\Lambda^q$
then there are infinitely many non-isomorphic simple representations of the same dimension as $X$.
Namely, the set $\mathcal{S}$ of simple representations forms an open subset of $\Rep(\Lambda^q,\alpha)$,
and since any simple representation has trivial endomorphism algebra, by Theorem~\ref{t:var},
$\mathcal{S}$ is either empty or smooth of dimension $g+2p(\alpha)$.
However, any orbit in $\mathcal{S}$ has dimension $g$.

For simplicity we now assume that the field $K$ has characteristic zero.
In Section~\ref{s:reptype} we combine Theorems~\ref{t:homform} and \ref{t:var}
with methods already developed in \cite{CBnorm} to study preprojective algebras,
and prove the following result.

\begin{thm}
\label{t:simpcrit}
Suppose that $\alpha$ and $q$ have the property that $p(\alpha) > p(\beta)+p(\gamma)+\dots$
for any nontrivial decomposition of $\alpha$ as a sum of positive
roots $\alpha = \beta+\gamma+\dots$ with $q^\beta = q^\gamma = \dots = 1$.
Then if $\Rep(\Lambda^q,\alpha)$ is non-empty, it is a complete intersection,
equidimensional of dimension $g+2p(\alpha)$, and the set of simple representations
is a dense open subset.
\end{thm}

We use this result in Section \ref{s:DSP} to prove Theorem \ref{t:DSP}.

Part of the writing up of this paper was done while the first author was visiting
the Mittag-Leffler Institute. He would like to thank his hosts for their hospitality.

\section{Independence of orientation and ordering}
\label{s:indep}
For $a\in \overline{Q}$, we define $g_a = 1 + aa^* \in K\overline{Q}$.
Note the following obvious formulas:
\begin{equation}
\label{e:swapg}
g_a a = a g_{a^*},
\quad
a^* g_a = g_{a^*} a^*.
\end{equation}
Let $L_Q$ be the algebra obtained from $K\overline{Q}$ by adjoining inverses
for the elements $g_a$.
(This is a trivial example of a universal localization, see for example \cite[Chapter 4]{Schofbook}.)

We identify $q\in (K^*)^I$ with the element $q = \sum_{v\in I} q_v e_v$ in $L_Q$.
Clearly it is invertible, with inverse $q^{-1} = \sum_{v\in I} q_v^{-1} e_v$.
We say that $x\in L_Q$ has \emph{diagonal Peirce decomposition} if
$x\in \bigoplus_{v\in I} e_vL_Q e_v$, or equivalently $x = \sum_{v\in I} e_v x e_v$.
For example $g_a$ has diagonal Peirce decomposition, and hence so also does $g_a^{-1}$.
Clearly $q$ commutes with any element which has diagonal Peirce decomposition.
Note also that if $x$ has diagonal Peirce decomposition and $y$ belongs to
$e_v L_Q e_w$, then so do $xy$ and $yx$.

The algebra $\Lambda^q(Q,<)$ is the quotient of $L_Q$ by the ideal generated by the element
\begin{equation}
\label{e:Lambdarelns}
\rho_{Q,<} = g_{a_1}^{\epsilon(a_1)} g_{a_2}^{\epsilon(a_2)} \dots g_{a_n}^{\epsilon(a_n)} - q,
\end{equation}
where $a_1 < a_2 < \dots < a_n$ are the arrows in $\overline{Q}$.

We prove Theorem~\ref{t:indeporor}. First independence of orientation.
Suppose that $a$ is an arrow in $Q$, and let
$Q'$ be the quiver obtained from $Q$ by deleting the arrow $a$, and
replacing it with a reverse arrow $b$, so $h(b)=t(a)$ and $t(b)=h(a)$.

There is an algebra homomorphism
$\theta:K\overline{Q'} \to L_{Q}$
sending the trivial paths $e_v$ and arrows other than $b$ and $b^*$ to themselves,
and sending $b$ to $a^*$ and $b^*$ to $-g_a^{-1} a$.
Clearly
\[
g_a \theta(g_{b^*}) = g_a (1 - g_{a}^{-1}aa^*) = g_a - aa^* = 1.
\]
Now by construction $g_a$ is invertible in $L_Q$, and hence
$g_a$ and $\theta(g_{b^*})$ are inverses in $L_Q$.
Also
\[
g_{a^*}\theta(g_b) = g_{a^*} (1-a^* g_a^{-1}a) = g_{a^*} - a^*a = 1
\]
using (\ref{e:swapg}). Thus $g_{a^*}$ and $\theta(g_b)$ are inverses in $L_Q$.
It follows that $\theta$ extends uniquely to a homomorphism
$\tilde{\theta}:L_{Q'}\to L_Q$.
To show it is an isomorphism we construct its inverse. As above, there is an algebra homomorphism
$\phi:K\overline{Q}\to L_{Q'}$
sending $a$ to $-b^*g_b^{-1}$ and $a^*$ to $b$.
One checks that $g_{b^*}$ and $\phi(g_a)$ are inverses
in $L_{Q'}$ and that $g_b$ and $\phi(g_{a^*})$ are inverses in
$L_{Q'}$. Thus $\phi$ extends to a homomorphism $\tilde{\phi}:L_Q\to L_{Q'}$.
Now $\tilde{\phi}(\tilde{\theta}(b)) = b$ and
\[
\tilde{\phi}(\tilde{\theta}(b^*)) = \tilde{\phi}(-g_a^{-1} a) = g_{b^*} b^* g_b^{-1} = b^*,
\]
while $\tilde{\theta}(\tilde{\phi}(a^*)) = a^*$ and
\[
\tilde{\theta}(\tilde{\phi}(a)) = \tilde{\theta}(-b^* g_b^{-1}) = g_a^{-1}a g_{a^*} = a,
\]
so that $\tilde{\phi}$ is the inverse to $\tilde{\theta}$.

Given an ordering $<$ on the arrows in $\overline{Q}$, let $<'$ be the
corresponding ordering for $\overline{Q'}$, with $a$ replaced by $b^*$ and $a^*$ replaced by $b$.
Clearly $\tilde{\theta}(\rho_{Q',<'}) = \rho_{Q,<}$, so $\Lambda^q(Q',<')\cong \Lambda^q (Q,<)$,
as required.

Now we show independence of the ordering.
Let $a_1 < \dots < a_n$ be the arrows in $\overline{Q}$.
Because $q$ commutes with any $g_a$, the relation $\rho_{Q,<}$
can be conjugated by $g_{a_1}^{\epsilon(a_1)}$ to give
\[
g_{a_2}^{\epsilon(a_2)} \dots g_{a_n}^{\epsilon(a_n)} g_{a_1}^{\epsilon(a_1)} - q,
\]
which shows that the algebra $\Lambda^q$ only depends on the induced \emph{cyclic} ordering
of $\overline{Q}$.

It thus suffices to show that $\Lambda^q(Q,<)\cong \Lambda^q(Q,<'')$, where
$<''$ is the ordering with the first two arrows exchanged, $a_2 <'' a_1 <'' a_3 <'' \dots <'' a_n$.
We may suppose that $h(a_1)= h(a_2)$, for otherwise $g_{a_1}$ and $g_{a_2}$ commute,
and $\rho_{Q,<} = \rho_{Q,<''}$.

If $a_1=a_2^*$, then $a_1$ is a loop, and $\Lambda^q(Q,<'')$ is the same as
the algebra \mbox{$\Lambda^q(Q',<')$} obtained by reversing the arrow $a_1$. The
argument above shows that this is isomorphic to $\Lambda^q(Q,<)$.

Thus suppose that $a_1\neq a_2^*$. By reversing arrows if necessary
we may assume that $\epsilon(a_1)=\epsilon(a_2)=1$. Define a homomorphism
$\theta:K\overline{Q}\to L_Q$
sending the trivial paths $e_v$ and arrows other than $a_1,a_1^*$ to themselves,
and with
\[
\theta(a_1)=g_{a_2} a_1,
\quad
\theta(a_1^*)=a_1^* g_{a_2}^{-1}.
\]
Clearly $\theta(g_b)=g_b$ for any arrow $b\neq a_1,a_1^*$. Moreover
\[
\theta(g_{a_1}) = g_{a_2} g_{a_1} g_{a_2}^{-1},
\quad
\theta(g_{a_1^*}) = g_{a_1^*}.
\]
Thus $\theta$ lifts to a homomorphism $\tilde{\theta}:L_Q\to L_Q$. Clearly this
is an isomorphism, and
\[
\tilde{\theta}(\rho_{Q,<}) =
\tilde{\theta}(g_{a_1} g_{a_2} g_{a_3}^{\epsilon(a_3)}\dots
g_{a_n}^{\epsilon(a_n)}  - q)
= g_{a_2} g_{a_1} g_{a_3}^{\epsilon(a_3)}\dots
g_{a_n}^{\epsilon(a_n)}  - q = \rho_{Q,<''},
\]
so that $\Lambda^q(Q,<'')\cong \Lambda^q(Q,<)$.

\section{Homological algebra}
\label{s:homform}
We retain the setup of Section~\ref{s:indep}, but simplify the notation,
writing $\Lambda$ rather than $\Lambda^q$ and $L$ instead of $L_Q$.
Thus $\Lambda = L/J$, where $J$ is the ideal generated by
\[
\rho = \biggl( \prod_{a\in \overline{Q}} g_a^{\epsilon(a)} \biggr) - q,
\]
or equivalently, since $\rho$ has diagonal Peirce decomposition,
by the elements $\rho_v = e_v \rho = \rho e_v$ ($v\in I$).

If $M$ is a $\Lambda$-$\Lambda$-bimodule, we say that elements $m_1,\dots,m_n\in M$ are
an \emph{$e$-free basis} of $M$ if each $m_i \in e_{v_i} M e_{w_i}$
for some vertices $v_i,w_i\in I$, and the natural homomorphism
\[
\bigoplus_{i=1}^n \Lambda e_{v_i} \otimes e_{w_i}\Lambda\to M
\]
sending $e_{v_i} \otimes e_{w_i}$ to $m_i$ is an isomorphism.
(Unadorned tensor products are over $K$.)
Let $P_0$ be a bimodule with $e$-free basis $\{ \eta_v \mid v\in I \}$,
where $\eta_v\in e_v P_0 e_v$ for all $v$, and let $P_1$ be a bimodule with
$e$-free basis $\{ \eta_a \mid a\in\overline{Q} \}$,
where $\eta_a\in e_{h(a)} P_1 e_{t(a)}$.

For $a\in\overline{Q}$ we define
\[
\ell_a = \prod_{\substack{b\in \overline{Q} \\ b<a}}
g_b^{\epsilon(b)},
\quad
r_a = \prod_{\substack{b\in \overline{Q} \\ b>a}}
g_b^{\epsilon(b)},
\]
both products taken in the appropriate order. The relation (\ref{c:rel}) can be written as
\begin{equation}
\label{e:ellgrisq}
\ell_a g_a^{\epsilon(a)} r_a = q,
\end{equation}
so that $r_a \ell_a = q g_a^{-\epsilon(a)}$.

Let $S$ be the semisimple subalgebra of $K\overline{Q}$ spanned by the trivial paths.
Clearly $L$ and $\Lambda$ are naturally $S$-rings.
Recall that if $A$ is an $S$-ring then there is a universal bimodule of derivations
$\Omega_S(A)$ which can be defined to be the kernel of the multiplication map $A\otimes_S A\to A$,
and the universal derivation $\delta_{A/S}:A\to \Omega_S(A)$
is then given by $\delta_{A/S}(a)=a\otimes 1-1\otimes a$. See \cite[Chapter 10]{Schofbook}.

\begin{lem}
\label{l:startprojres}
There is an exact sequence of $\Lambda$-$\Lambda$-bimodules
\[
P_0
\xrightarrow{\alpha}
P_1
\xrightarrow{\beta}
P_0
\xrightarrow{\gamma}
\Lambda \to 0
\]
where $\gamma(\eta_v) = e_v$,
$\beta(\eta_a) = a\eta_{t(a)} - \eta_{h(a)} a$, and
\[
\alpha(\eta_v) =\sum_{\substack{a\in \overline{Q} \\ h(a)=v}} \ell_a \Delta_a r_a
\]
where
\[
\Delta_a = \begin{cases}
\eta_a a^* + a \eta_{a^*}
& \text{\textup{(if $\epsilon(a)=1$)}} \\
- g_a^{-1} (\eta_a a^* + a \eta_{a^*}) g_a^{-1}
& \text{\textup{(if $\epsilon(a)=-1$).}}
\end{cases}
\]
\end{lem}

\begin{proof}
We consider the diagram of $\Lambda$-$\Lambda$-bimodules and homomorphism
\[
\begin{CD}
P_0 @>\alpha>> P_1 @>\beta>> P_0 @>\gamma>> \Lambda @>>> 0
\\
@V\theta VV @V\phi VV @V\psi VV @|
\\
J/J^2 @>\sigma>> \Lambda\otimes_L \Omega_S(L)\otimes_L \Lambda @>\zeta>> \Lambda\otimes_S \Lambda @>\mu>> \Lambda @>>> 0.
\end{CD}
\]
The second row is the exact sequence given by splicing the sequence of \cite[Theorem 10.3]{Schofbook}
with the defining sequence for $\Omega_S(\Lambda)$. Thus $\mu$ is the multiplication map, $\zeta$ is the natural map
\[
\Lambda\otimes_L \Omega_S(L)\otimes_L \Lambda \to
\Lambda\otimes_L (L\otimes_S L)\otimes_L \Lambda \cong
\Lambda\otimes_S \Lambda,
\]
and $\sigma(J^2 + x) = 1\otimes\delta_{L/S}(x)\otimes 1$ for $x\in J$.

Let $\psi$ be the isomorphism sending $\eta_v$ to $e_v\otimes e_v = e_v\otimes 1 = 1\otimes e_v$.

Let $B$ be the $S$-$S$-sub-bimodule of $K\overline{Q}$ spanned by the arrows,
so that $K\overline{Q}$ is identified with the tensor algebra of $B$ over $S$.
By \cite[Theorem 10.5]{Schofbook}, there is an isomorphism
\[
\Omega_S(K\overline{Q}) \cong K\overline{Q} \otimes_S B \otimes_S K\overline{Q}
\]
under which $\delta_{K\overline{Q}/S}(a)$ corresponds to $1\otimes a\otimes 1$ for $a\in\overline{Q}$.
Since $L$ is a universal localization of $K\overline{Q}$,
by \cite[Theorem 10.6]{Schofbook}, there is an isomorphism
\[
\Omega_S(L) \cong L \otimes_{K\overline{Q}} \Omega_S(K\overline{Q}) \otimes_{K\overline{Q}} L
\]
under which $\delta_{L/S}(a)$ corresponds to $1\otimes\delta_{K\overline{Q}/S}(a)\otimes 1$ for $a\in\overline{Q}$.
These give an isomorphism $\phi$,
\[
P_1 \cong \Lambda \otimes_S B \otimes_S \Lambda \cong \Lambda\otimes_L \Omega_S(L) \otimes_L \Lambda.
\]
Thus $\phi(\eta_a) = 1\otimes\delta(a)\otimes 1$
for $a\in\overline{Q}$.

Let $\theta$ be the homomorphism sending $\eta_v$ to $\rho_v$. Since $J$ is generated by the elements
$\rho_v$, it follows that $\theta$ is onto.

We check that the diagram commutes. Only the left-hand square is non-trivial.
Since $\delta_{L/S}$ is a derivation,
\[
\delta_{L/S}(g_a) = \delta_{L/S}(1+aa^*) = \delta_{L/S}(a)a^* + a\delta_{L/S}(a^*)
\]
and, by considering $\delta_{L/S}(g_a g_a^{-1})$,
\[
\delta_{L/S}(g_a^{-1}) = - g_a^{-1} \left( \delta_{L/S}(a)a^* + a\delta_{L/S}(a^*) \right) g_a^{-1}.
\]
Thus $1\otimes\delta_{L/S}(g_a^{\epsilon(a)})\otimes 1 = \phi(\Delta_a)$.
Also
\[
\delta_{L/S} (\rho) = \delta_{L/S} \left(\prod_{a\in \overline{Q}} g_a^{\epsilon(a)} - q \right)
= \sum_{a\in \overline{Q}} \ell_a \delta_{L/S} (g_a^{\epsilon(a)}) r_a.
\]
Since $e_v$ commutes with the elements $g_a$, and $\delta_{L/S}(e_v)=0$, we deduce that
\[
\delta_{L/S} (\rho_v) = \sum_{\substack{a\in \overline{Q} \\ h(a)=v}} \ell_a \delta_{L/S} (g_a^{\epsilon(a)}) r_a.
\]
Thus
\[
1\otimes\delta_{L/S} (\rho_v)\otimes 1 = \sum_{\substack{a\in \overline{Q} \\ h(a)=v}} \ell_a \phi(\Delta_a) r_a = \phi(\alpha(\eta_v)),
\]
so that $\phi\alpha=\sigma\theta$, as required.
It follows that the top row is exact.
\end{proof}

We consider $\Lambda\otimes\Lambda$ as a $\Lambda$-$\Lambda$-bimodule, with the action given by
\[
\lambda (x\otimes x')\lambda' = \lambda x\otimes x'\lambda'
\]
for $\lambda,\lambda'\in\Lambda$ and $x\otimes x'\in \Lambda\otimes\Lambda$.
There is a duality $P\rightsquigarrow P^\vee = \Hom_{\Lambda-\Lambda}(P,\Lambda\otimes\Lambda)$
on the category of finitely generated projective $\Lambda$-$\Lambda$-bimodules,
where the bimodule structure on $P^\vee$ is given by
\begin{equation}
\label{e:dualbmodstr}
(\lambda f \lambda')(p) = \sum_j x_j\lambda' \otimes \lambda x'_j
\end{equation}
for $\lambda,\lambda'\in\Lambda$, $f\in P^\vee$ and $p\in P$, where $f(p) = \sum_j x_j\otimes x'_j$.
Observe that
\begin{equation}
\label{e:bimoddual}
(\Lambda e_v \otimes e_w\Lambda)^\vee \cong \Lambda e_w\otimes e_v \Lambda,
\end{equation}
with $a\otimes b\in \Lambda e_w\otimes e_v \Lambda$
corresponding to the homomorphism sending $e_v\otimes e_w$ to $b\otimes a$.
Thus the dual of an $e$-free bimodule is $e$-free. If $P$ has $e$-free basis $m_1,\dots,m_n$
with $m_i\in e_{v_i} P e_{w_i}$, then $P^\vee$ has $e$-free basis
$m_1^\vee,\dots,m_n^\vee$ defined by
\[
m_i^\vee(m_j) = \begin{cases}
e_{v_i}\otimes e_{w_i} &\text{(if $i=j$)}
\\
0 &\text{(if $i\neq j$).}
\end{cases}
\]
In view of (\ref{e:dualbmodstr}), one has $m_i^\vee \in e_{w_i} P^\vee e_{v_i}$.

For $a\in\overline{Q}$ we define elements $c_a = q^{-1} \ell_a a r_{a^*}$ in $\Lambda$,
\[
\xi_a = \begin{cases}
q^{-1} \ell_{a^*} \eta_a^\vee g_a r_a
&
\text{(if $\epsilon(a)=1$)}
\\
-q^{-1} \ell_{a^*} g_{a^*} \eta_a^\vee r_a
&
\text{(if $\epsilon(a)=-1$).}
\end{cases}
\]
in $P_1^\vee$ and $\xi_v = q \eta_v^\vee = \eta_v^\vee q$ in $P_0^\vee$.

\begin{lem}
\label{l:cechcalc}
We have $\alpha^\vee(\xi_a) = c_{a^*} \xi_{h(a)}- \xi_{t(a)} c_{a^*}$.
\end{lem}

\begin{proof}
From the definition of $\alpha$ we have
\[
\begin{split}
\alpha(\eta_v) =
&\sum_{\substack{h(a)=v \\ \epsilon(a)=1 }} \ell_a \eta_a a^* r_a
+ \sum_{\substack{t(a)=v \\ \epsilon(a)=-1 }} \ell_{a^* }a^* \eta_a r_{a^*}
\\
&- \sum_{\substack{h(a)=v \\ \epsilon(a)=-1 }} \ell_a g_a^{-1} \eta_a a^* g_a^{-1} r_a
- \sum_{\substack{t(a)=v \\ \epsilon(a)=1 }} \ell_{a^*} g_{a^*}^{-1} a^* \eta_a g_{a^*}^{-1} r_{a^*}.
\end{split}
\]
Extracting the terms involving $\eta_a$ in this expression,
we get
\[
\alpha^\vee(\eta_a^\vee) =
\begin{cases}
a^* r_a \eta_{h(a)}^\vee \ell_a - g_{a^*}^{-1} r_{a^*} \eta_{t(a)}^\vee \ell_{a^*} g_{a^*}^{-1} a^*
&
\text{(if $\epsilon(a)=1$)}
\\
r_{a^*} \eta_{t(a)}^\vee \ell_{a^*} a^* - a^* g_a^{-1} r_a \eta_{h(a)}^\vee \ell_a g_a^{-1}
&
\text{(if $\epsilon(a)=-1$).}
\end{cases}
\]
Then if $\epsilon(a)=1$ we have
\begin{align*}
\alpha^\vee(\xi_a)
&= \alpha^\vee(q^{-1} \ell_{a^*} \eta_a^\vee g_a r_a)
\\
&= q^{-1} \ell_{a^*} \alpha^\vee(\eta_a^\vee) g_a r_a
\\
&= q^{-1} \ell_{a^*} (a^* r_a \eta_{h(a)}^\vee \ell_a - g_{a^*}^{-1} r_{a^*} \eta_{t(a)}^\vee \ell_{a^*} g_{a^*}^{-1} a^*) g_a r_a
\\
&= q^{-1} \ell_{a^*} a^* r_a \eta_{h(a)}^\vee q - \eta_{t(a)}^\vee \ell_{a^*} g_{a^*}^{-1} a^* g_a r_a
&\text{by (\ref{e:ellgrisq}),}
\\
&= q^{-1} \ell_{a^*} a^* r_a \eta_{h(a)}^\vee q - \eta_{t(a)}^\vee \ell_{a^*} a^* r_a
&\text{by (\ref{e:swapg}),}
\\
&= c_{a^*} \xi_{h(a)} - \xi_{t(a)} c_{a^*}.
\end{align*}
A similar calculation gives the result if $\epsilon(a)=-1$.
\end{proof}

\begin{lem}
\label{l:funnyonto}
The algebra homomorphism $\theta:K\overline{Q}\to \Lambda$ defined by $\theta(e_v)=e_v$ for $v\in I$
and $\theta(a) = c_a$ for $a\in\overline{Q}$
induces a surjective homomorphism $\tilde{\theta}:L\to \Lambda$.
\end{lem}

\begin{proof}
We have
\begin{align*}
\theta(g_a)
&= 1 + c_a c_{a^*}
\\
&= 1 + q^{-1} \ell_a a r_{a^*} q^{-1} \ell_{a^*} a^* r_a
\\
&= 1 + q^{-1} \ell_a a g_{a^*}^{-\epsilon(a^*)} a^* r_a
&\text{by (\ref{e:ellgrisq}),}
\\
&= 1 + q^{-1} \ell_a a a^* g_{a}^{\epsilon(a)} r_a
&\text{by (\ref{e:swapg}),}
\\
&= 1 + q^{-1} \ell_a a a^* \ell_a^{-1} q
&\text{by (\ref{e:ellgrisq}),}
\\
&= 1 + \ell_a a a^* \ell_a^{-1}
\\
&= \ell_a g_a \ell_a^{-1}.
\end{align*}
Since this is invertible, there is an induced homomorphism $\tilde{\theta}:L\to \Lambda$.
The image of $\tilde{\theta}$ contains $\ell_a g_a \ell_a^{-1}$
and its inverse $\ell_a g_a^{-1} \ell_a^{-1}$.
For $a$ minimal with respect to the ordering,
we have $\ell_a=1$, so these elements are $g_a$ and $g_a^{-1}$. Then, by induction
working up the ordering, the image contains $g_a^{\pm1}$ for all $a$.
Thus the image contains $\ell_a^{\pm1}$ and $r_a^{\pm1}$.
Then, since the image contains $\ell_a a r_{a^*}$, it contains $a$. Thus $\tilde{\theta}$ is onto.
\end{proof}

\begin{lem}
\label{l:cexactseq}
There is an exact sequence
\[
P_1 \xrightarrow{\phi} P_0 \xrightarrow{\psi} \Lambda \to 0
\]
where $\psi(\eta_v) = e_v$ and $\phi(\eta_a) = c_a \eta_{t(a)}- \eta_{h(a)}c_a$.
\end{lem}

\begin{proof}
We consider $\Lambda$ as a left or right $L$-module using
the surjective homomorphism of Lemma~\ref{l:funnyonto}.
Inducing up the defining sequence for $\Omega_S(L)$ gives an
exact sequence
\[
\Lambda\otimes_L \Omega_S(L)\otimes_L \Lambda \xrightarrow{\omega}
\Lambda\otimes_S \Lambda \xrightarrow{\mu} \Lambda \to 0
\]
where $\mu$ is multiplication. The map $\omega$ sends $1\otimes(\sum_k x_k\otimes x_k')\otimes 1$
to $\sum_k \tilde{\theta}(x_k)\otimes \tilde{\theta}(x_k)$.
Now $\Lambda\otimes_S \Lambda \cong P_0$, with $\eta_v$ corresponding to $e_v\otimes e_v$, and
$\Lambda\otimes_L \Omega_S(L)\otimes_L \Lambda\cong P_1$
with $\eta_a$ corresponding to $1\otimes \delta_{L/S}(a) \otimes 1$.
Now
\begin{align*}
\omega(1\otimes\delta_{L/S}(a)\otimes 1)
&= \omega(1\otimes(a\otimes 1-1\otimes a)\otimes 1)
\\
&= \tilde{\theta}(a)\otimes 1 - 1\otimes \tilde{\theta}(a)
\\
&= c_a \otimes 1 - 1\otimes c_a
\\
&= c_a (e_{t(a)}\otimes e_{t(a)}) - (e_{h(a)}\otimes e_{h(a)}) c_a
\end{align*}
The lemma follows.
\end{proof}

\begin{lem}
\label{l:magicseq}
%The cokernel of $\alpha^\vee$ is isomorphic to $\Lambda$ as a $\Lambda$-$\Lambda$-bimodule.
There is an exact sequence $P_1^\vee\xrightarrow{\alpha^\vee}P_0^\vee\to\Lambda\to 0$.
\end{lem}

\begin{proof}
Observe that the $\xi_v$ ($v\in I$) are an $e$-free basis of $P_0^\vee$ and
the $\xi_a$ ($a\in\overline{Q}$) are an $e$-free basis of $P_1^\vee$.
Let $f$ be the isomorphism $P_0\to P_0^\vee$ sending $\eta_v$ to
$\xi_v$, and let $g$ be the isomorphism $P_1\to P_1^\vee$ sending $\eta_a$ to $\xi_{a^*}$.
In view of Lemma~\ref{l:cexactseq}, it suffices to prove that $\alpha^\vee g = f\phi$.
This follows from Lemma~\ref{l:cechcalc}.
\end{proof}

We now turn to the proof of Theorem~\ref{t:homform}, which is analogous to
that of \cite[Lemma 1]{CBex}. From Lemma~\ref{l:startprojres} we have the start of a projective resolution of $X$,
\[
P_0\otimes_\Lambda X \to P_1\otimes_\Lambda X \to P_0\otimes_\Lambda X \to X\to 0.
\]
Applying $\Hom_\Lambda(-,Y)$ gives a complex
\begin{equation}
\label{e:LambdaXYcplx}
0 \to \Hom_\Lambda(P_0\otimes_\Lambda X,Y) \to \Hom_\Lambda(P_1\otimes_\Lambda X,Y)\to \Hom_\Lambda(P_0\otimes_\Lambda X,Y) \to 0
\end{equation}
such that the cohomology at the first two places is $\Hom(X,Y)$ and $\Ext^1(X,Y)$.
To understand the cohomology at at the third place we dualize to give
\[
0 \to \Hom_\Lambda(P_0\otimes_\Lambda X,Y)^* \to \Hom_\Lambda(P_1\otimes_\Lambda X,Y)^* \to \dots
\]
For $P$ a finitely generated projective $\Lambda$-$\Lambda$-bimodule there is a natural isomorphism
\[
\Hom_\Lambda(P\otimes_\Lambda Y,X)
\cong \Hom_\Lambda(Y,\Hom_{\Lambda-\Lambda}(P,X))
\cong \Hom_\Lambda(Y,P^\vee\otimes_\Lambda X)
\]
This gives a natural transformation
\[
\Hom_\Lambda(P\otimes_\Lambda Y,X) \to \Hom_\Lambda(P^\vee\otimes_\Lambda X,Y)^*,
\quad
f \mapsto (g\mapsto \tr(g f'))
\]
where $f'\in \Hom_\Lambda(Y,P^\vee\otimes_\Lambda X)$ corresponds to $f$.
Clearly this is an isomorphism in case $P=\Lambda\otimes\Lambda$, so
it is a natural isomorphism for all finitely generated projective $P$.
Using this we rewrite the dualized complex as
\[
0 \to \Hom_\Lambda(P_0^\vee\otimes_\Lambda Y,X) \to \Hom_\Lambda(P_1^\vee\otimes_\Lambda Y,X) \to \dots.
\]
Using Lemma~\ref{l:magicseq} we see that the cohomology in the first position is $\Hom_\Lambda(Y,X)$.
Now the alternating sum of the dimensions of the cohomology spaces in (\ref{e:LambdaXYcplx}),
\[
\dim\Hom_\Lambda(X,Y) - \dim\Ext^1_\Lambda(X,Y) + \dim \Hom_\Lambda(Y,X),
\]
is equal to
\[
2\dim \Hom(P_0\otimes_\Lambda X,Y) - \dim \Hom(P_1\otimes_\Lambda X,Y),
\]
which is $(\underline{\dim}X,\underline{\dim}Y)$.

\section{Middle convolution}
\label{s:MC}
In this section we prove Theorem~\ref{t:MC}.
Let $v$ be a loopfree vertex in $Q$ and suppose that $q_v\neq 1$.
Note that by reorienting we may assume that no arrow in $Q$ has tail at $v$.
Suppose that the arrows with head at $v$ are $a_1<a_2<\dots<a_n$.
Let $q'=u_v(q)$.

Given a representation $X$ of $\Lambda^q$, we consider it as a representation of $\overline{Q}$
satisfying (\ref{c:invtwoprime}) and (\ref{c:reltwoprime}).
For $1\le i\le n+1$ define
\[
\ell_i = (1_{X_v}+X_{a_1}X_{a_1^*})(1_{X_v}+X_{a_2}X_{a_2^*}) \dots (1_{X_v}+X_{a_{i-1}}X_{a_{i-1}^*}).
\]
Clearly
\begin{equation}
\label{e:lXXeq}
\sum_{j=1}^{i-1} \ell_j X_{a_j} X_{a_j^*} = \ell_i - 1_{X_v}
\end{equation}
and using the relation (\ref{c:reltwoprime}),
\begin{equation}
\label{e:lnqm}
\sum_{j=1}^{n} \ell_j X_{a_j} X_{a_j^*} = \ell_{n+1} - 1_{X_v} = (q_v-1) 1_{X_v}
\end{equation}
Define
\[
X_\oplus = \bigoplus_{i=1}^n X_{t(a_i)}.
\]
Let $\iota_i:X_{t(a_i)}\to X_\oplus$ and $\pi_i:X_\oplus\to X_{t(a_i)}$ be the natural
maps, and define
\begin{equation}
\label{e:iotapidef}
\iota = \sum_{i=1}^n \iota_i X_{a_i^*} : X_v\to X_\oplus,
\quad
\pi = \frac{1}{q_v-1} \sum_{i=1}^n \ell_i X_{a_i} \pi_i:X_\oplus\to X_v,
\end{equation}
Equation (\ref{e:lnqm}) ensures that $\pi\iota=1_{X_v}$.
Thus $\iota\pi$ and $\epsilon = 1_{X_\oplus}-\iota\pi$ are idempotent endomorphisms of $X_\oplus$.
Define
\[
\phi_i : X_{t(a_i)} \to X_\oplus,
\quad
\phi_i = \sum_{j=1}^{i-1} \iota_j X_{a_j^*}X_{a_i}
+ \frac{1}{q_v} \sum_{j=i}^n \iota_j X_{a_j^*} X_{a_i}
+ \frac{1-q_v}{q_v} \iota_i.
\]
Observe that
\begin{equation}
\label{e:piphiless}
\pi_j \phi_i = X_{a_j^*}X_{a_i} \quad\text{for $j<i$.}
\end{equation}

\begin{lem}
\label{l:piphiz}
$\pi \phi_i=0$.
\end{lem}

\begin{proof}
We have
\begin{align*}
\pi \phi_i &= \frac{1}{q_v-1} \sum_{k=1}^n \ell_k X_{a_k} \pi_k \phi_i
\\
&= \frac{1}{q_v-1} \sum_{k=1}^n \ell_k X_{a_k} \pi_k
\left(
\sum_{j=1}^{i-1} \iota_j X_{a_j^*}X_{a_i}
+ \frac{1-q_v}{q_v} \iota_i
+ \frac{1}{q_v} \sum_{j=i}^n \iota_j X_{a_j^*} X_{a_i}
\right)
\\
&= \frac{1}{q_v-1}
\left(
\sum_{j=1}^{i-1} \ell_j X_{a_j} X_{a_j^*}X_{a_i}
+ \frac{1-q_v}{q_v} \ell_i X_{a_i}
+ \frac{1}{q_v} \sum_{j=i}^n \ell_j X_{a_j} X_{a_j^*} X_{a_i}
\right)
\\
&= \frac{1}{q_v-1}
\left(
\sum_{j=1}^{i-1} \ell_j X_{a_j} X_{a_j^*}
+ \frac{1-q_v}{q_v} \ell_i
+ \frac{1}{q_v} \sum_{j=i}^n \ell_j X_{a_j} X_{a_j^*}
\right) X_{a_i}.
\end{align*}
Simplifying this using (\ref{e:lXXeq}) and
\[
\sum_{j=i}^{n} \ell_j X_{a_j} X_{a_j^*}
= (q_v-1)1_{X_v} - \sum_{j=1}^{i-1} \ell_j X_{a_j} X_{a_j^*}
= q_v 1_{X_v} - \ell_i,
\]
the lemma follows.
\end{proof}

\begin{lem}
\label{l:phipiprod}
For all $0\le m\le n$ we have
\[
(1_{X_\oplus} + \phi_1\pi_1)(1_{X_\oplus} + \phi_2\pi_2)\dots(1_{X_\oplus} + \phi_m\pi_m)
=
1_{X_\oplus} + \frac{1-q_v}{q_v} \sum_{j=1}^{m} \epsilon \iota_j \pi_j.
\]
\end{lem}

\begin{proof}
We prove the assertion by induction on $m$. It is trivially true for $m=0$,
and assuming the truth for $m-1$, to deduce it for $m$ we need to show that
\[
\left(
1_{X_\oplus} + \frac{1-q_v}{q_v} \sum_{j=1}^{m-1} \epsilon \iota_j \pi_j
\right)
(1_{X_\oplus} + \phi_m\pi_m)
=
1_{X_\oplus} + \frac{1-q_v}{q_v} \sum_{j=1}^{m} \epsilon \iota_j \pi_j,
\]
or equivalently that
\begin{equation}
\label{e:phipineed}
\phi_m\pi_m = \frac{1-q_v}{q_v} \epsilon
\left(
\iota_m \pi_m - \sum_{j=1}^{m-1} \iota_j \pi_j \phi_m \pi_m
\right).
\end{equation}
Now by (\ref{e:piphiless}), the right hand side of this is
\[
\frac{1-q_v}{q_v} (1-\iota\pi)
\left(
\iota_m \pi_m - \sum_{j=1}^{m-1} \iota_j X_{a_j^*} X_{a_m} \pi_m
\right).
\]
Multiplying out and using that $\pi\iota_j = \frac{1}{q_v-1} \ell_j X_{a_j}$ this gives
\[
\frac{1-q_v}{q_v} \left(
\iota_m \pi_m - \sum_{j=1}^{m-1} \iota_j X_{a_j^*} X_{a_m} \pi_m
\right)
+ \frac{1}{q_v} \left(
\iota \ell_m X_{a_m} \pi_m - \sum_{j=1}^{m-1} \iota \ell_j X_{a_j} X_{a_j^*} X_{a_m} \pi_m
\right).
\]
Thanks to (\ref{e:lXXeq}), this becomes
\[
\frac{1-q_v}{q_v} \left(
\iota_m \pi_m - \sum_{j=1}^{m-1} \iota_j X_{a_j^*} X_{a_m} \pi_m
\right)
+ \frac{1}{q_v} \iota X_{a_m} \pi_m.
\]
Now expanding using the formula for $\iota$ and rearranging, this gives $\phi_m \pi_m$, proving
(\ref{e:phipineed}), as required.
\end{proof}

Let $X'$ be the representation of $\overline{Q}$ defined as follows.
The vector spaces are $X'_w = X_w$ for vertices $w\neq v$ and
$X'_v = \Ima(\epsilon) = \Ker(\iota\pi) = \Ker(\pi)$.
Let $\iota'$ be the inclusion $X'_v\to X_\oplus$. The maps are
$X'_a = X_a$ for arrows $a\in \overline{Q}$ not incident at $v$,
$X'_{a_i^*} = \pi_i \iota'$ and $X'_{a_i}: X'_{t(a_i)}\to X'_v$ is the map with $\phi_i=\iota' X'_{a_i}$.
It exists by Lemma~\ref{l:piphiz}, and is unique since $\iota'$ is injective.

\begin{lem}
$X'$ is a representation of $\Lambda^{q'}$. If $X$ has dimension $\alpha$, then $X'$ has dimension $s_v(\alpha)$.
\end{lem}

\begin{proof}
It is clear that
\[
1_{X_{t(a_i)}} + X'_{a_i^*} X'_{a_i} = \frac{1}{q_v} (1_{X_{t(a_i)}} + X_{a_i^*} X_{a_i}),
\]
which implies that $X'$ satisfies the relations (\ref{c:reltwoprime}) at vertices different from $v$.
On the other hand, taking $m=n$ in Lemma~\ref{l:phipiprod}, one has
\[
(1_{X_\oplus} + \phi_1\pi_1)\dots(1_{X_\oplus} + \phi_n\pi_n)
=
1_{X_\oplus} + \frac{1-q_v}{q_v} \epsilon,
\]
which on restricting to $X'_v$ gives
\[
(1_{X'_v} + X'_{a_1} X'_{a_1^*}) \dots (1_{X'_v} + X'_{a_n} X'_{a_n^*}) = \frac{1}{q_v} 1_{X'_v},
\]
which shows that (\ref{c:invtwoprime}) and the relation (\ref{c:reltwoprime}) holds at the vertex $v$.
Thus $X'$ is a representation of $\Lambda^{q'}$.
The assertion about dimension vectors is obvious, since $\dim X_v' = \dim X_\oplus - \dim X_v$.
\end{proof}

\begin{lem}
The assignment $X \rightsquigarrow X'$ defines a functor $F_q$ from representations of $\Lambda^q$
to representations of $\Lambda^{q'}$. It is an equivalence, with inverse $F_{q'}$.
\end{lem}

\begin{proof}
It is clear that the construction of $X'$ defines a functor $F_q$.
The analogous functor $F_{q'}$ applied to $X'$ defines a representation $X''$ of $\Lambda^q$.
We show that $X''$ is naturally isomorphic to $X$.

Note that $X'_\oplus = X_\oplus$.
Let $\ell_i'$ and $\pi'$ be the analogues of $\ell_i$ and $\pi$, but constructed from $X'$.
We have already defined a map $\iota'$, and this is the analogue of $\iota$ since
\[
\sum_{i=1}^n \iota_i X'_{a_i^*}
=
\sum_{i=1}^n \iota_i \pi_i \iota' = \iota'.
\]

We have
\begin{align*}
\iota'\pi'
&= \frac{1}{q'_v-1} \sum_{i=1}^n \iota' \ell'_i X'_{a_i} \pi_i
\\
&= \frac{q_v}{1-q_v} \sum_{i=1}^n \iota' (1_{X_\oplus}+\phi_1\pi_1) \dots (1_{X_\oplus}+\phi_{i-1}\pi_{i-1}) \phi_i \pi_i
\\
&= \frac{q_v}{1-q_v} \sum_{i=1}^n \phi_i \pi_i + \sum_{i=1}^n \sum_{j=1}^{i-1} \epsilon \iota_j \pi_j \phi_i \pi_i
\end{align*}
by Lemma~\ref{l:phipiprod}. Now using (\ref{e:piphiless}) and
\[
\epsilon \iota_j
= \iota_j - \iota\pi\iota_j
= \iota_j - \frac{1}{q_v-1} \iota \ell_j X_{a_j}
= \iota_j - \frac{1}{q_v-1} \sum_{k=1}^n \iota_k X_{a_k^*} \ell_j X_{a_j},
\]
we obtain
\[
\iota'\pi'
= \frac{q_v}{1-q_v} \sum_{i=1}^n \phi_i \pi_i
+ \sum_{i=1}^n \sum_{j=1}^{i-1} \iota_j X_{a_j^*} X_{a_i} \pi_i
- \frac{1}{q_v-1}\sum_{i=1}^n \sum_{j=1}^{i-1}\sum_{k=1}^n \iota_k X_{a_k^*} \ell_j X_{a_j} X_{a_j^*} X_{a_i} \pi_i.
\]
By (\ref{e:lXXeq}) this gives
\[
\iota'\pi'
= \frac{q_v}{1-q_v} \sum_{i=1}^n \phi_i \pi_i
+ \sum_{i=1}^n \sum_{j=1}^{i-1} \iota_j X_{a_j^*} X_{a_i} \pi_i
- \frac{1}{q_v-1}\sum_{i=1}^n \sum_{k=1}^n \iota_k X_{a_k^*} (\ell_i - 1) X_{a_i} \pi_i.
\]
Expanding this with the formula for $\phi_i$, most terms cancel, leaving
\[
\iota'\pi'
= 1_{X_\oplus} + \frac{q_v}{1-q_v} \sum_{i=1}^n \sum_{k=1}^n \iota_k X_{a_k^*} \ell_i X_{a_i} \pi_i = 1_{X_\oplus} - \iota\pi.
\]
Thus $X''_v$ can be naturally identified with $X_v$, with the inclusion $\iota''$ of $X''_v$
into $X_\oplus$ then identified with $\iota$.
The linear maps defining $X''$ are then given by
$X''_{a_i^*} = \pi_i \iota = X_{a_i^*}$ and
\begin{align*}
\iota X''_{a_i}
&= \sum_{j=1}^{i-1} \iota_j X'_{a_j^*} X'_{a_i} + \frac{1}{q'_v} \sum_{j=i}^{n} \iota_j X'_{a_j^*} X'_{a_i} + \frac{1-q'_v}{q'_v} \iota_i
\\
&= \sum_{j=1}^{i-1} \iota_j \pi_j \iota' X'_{a_i} + q_v \sum_{j=i}^{n} \iota_j \pi_j \iota' X'_{a_i} + (q_v-1) \iota_i
\\
&= \sum_{j=1}^{i-1} \iota_j \pi_j \phi_i + q_v \sum_{j=i}^{n} \iota_j \pi_j \phi_i + (q_v-1) \iota_i
\\
&= \sum_{j=1}^{i-1} \iota_j X_{a_j^*} X_{a_i}
+ q_v \left( \sum_{j=i}^{n} \iota_j \frac{1}{q_v} X_{a_j^*} X_{a_i} + \frac{1-q_v}{q_v} \iota_i \right)
+ (q_v-1) \iota_i
\\
&= \sum_{j=1}^n \iota_j X_{a_j^*} X_{a_i} = \iota X_{a_i}.
\end{align*}
Thus $X'' = X$, as desired.
\end{proof}

This completes the proof of Theorem~\ref{t:MC}.

\section{Simple representations}
\label{s:simp}
We begin with a lemma. Compare it with \cite[Theorem 1]{Scott} and \cite[Lemma 7.2]{CBmm}.

\begin{lem}
\label{l:simpcond}
If $\Lambda^q$ has a simple representation of dimension $\alpha$ and $v$ is a vertex, then
either $\alpha=\epsilon_v$ or $q_v\neq 1$ or $(\alpha,\epsilon_v)\le 0$.
\end{lem}

\begin{proof}
Suppose otherwise. Since $(\alpha,\epsilon_v)>0$ there is no loop at $v$.
Suppose that the arrows with head at $v$ are $a_1<a_2<\dots<a_n$.
By reorienting we may assume that $\epsilon(a_i)=1$ for all $i$.
Let $X$ be a simple representation of dimension $\alpha$.
Since $q_v=1$ the relation at vertex $v$ can be rewritten as
\[
\sum_{i=1}^n \ell_{a_i} X_{a_i} X_{a_i^*} = 0
\]
where $\ell_{a_i} = \prod_{j=1}^{i-1} (1_{X_v} + X_{a_j} X_{a_j^*})$.
Define
\[
X_\oplus = \bigoplus_{i=1}^n X_{t(a_i)},
\]
let $\theta:X_v\to X_\oplus$ be the linear map with components $X_{a_i^*}$, and
let $\phi:X_\oplus\to X_v$ be the linear map with components $\ell_{a_i} X_{a_i}$.
Thus $\phi\theta=0$.

Suppose that $\theta$ is not injective. Then $X$ has a subrepresentation given
by the vector space $\Ker(\theta)$ at vertex $v$ and the zero subspace at all
other vertices. By simplicity $X$ is equal to this subrepresentation. But,
since there is no loop at $v$, the fact that $X$ is simple implies that
its dimension vector is $\epsilon_v$, a contradiction. Thus $\theta$ is injective.

Suppose that $\phi$ is not surjective. Then $X$ has a subrepresentation given
by the vector space $U=\Ima(\phi)$ at vertex $v$ and the whole vector space $X_w$
at all other vertices $w$. Namely, it suffices to prove that $\Ima(X_{a_i})\subseteq U$
for all $i$. Now we know that $\Ima(\ell_{a_i}X_{a_i})\subseteq U$. In case $i=1$
we have $\ell_{a_1}=1$, so this already gives $\Ima(X_{a_1})\subseteq U$.
It follows that $(1+X_{a_1}X_{a_1}^*)(U)\subseteq U$, and since $1+X_{a_1}X_{a_1}^*$
acts invertibly on $X_v$, we have $(1+X_{a_1}X_{a_1}^*)^{-1}(U)=U$.
Now $\Ima(\ell_{a_2} X_{a_2})\subseteq U$, and hence
\[
\Ima(X_{a_2}) \subseteq \ell_{a_2}^{-1} (U) = (1+X_{a_1}X_{a_1}^*)^{-1}(U)=U.
\]
Repeating in this way, one has $\Ima(X_{a_i})\subseteq U$ for all $i$, as required.
As above, the simplicity assumption leads to a contradiction, so that $\phi$ is surjective.

It follows that $\phi$ induces a surjective linear map $X_\oplus/\Ima(\theta)\to X_v$,
so that $\dim X_\oplus \ge 2\dim X_v$, and hence $(\alpha,\epsilon_v)\le 0$.
\end{proof}

We prove Theorem~\ref{t:simps} by induction, supposing its truth for all $\beta<\alpha$.
If $\alpha$ is in the fundamental region then it is an imaginary root, as required.
Thus suppose that $\alpha$ is not in the fundamental region.
Clearly $\alpha$ has connected support since there is a simple representation of
dimension $\alpha$.
Thus $(\alpha,\epsilon_v)>0$ for some vertex $v$. Then $v$ must be loopfree and
we have $s_v(\alpha)<\alpha$.

If $q_v\neq 1$ then Theorem~\ref{t:MC} shows that representations of $\Lambda^q$ of
dimension $\alpha$ correspond to representations of $\Lambda^{u_v(q)}$ of dimension $s_v(\alpha)$.
In particular $X$ corresponds to a simple representation.
Thus by induction $s_v(\alpha)$ is a root, and hence so is $\alpha$.
Moreover, if $\alpha$ is a real root, then so is $s_v(\alpha)$, so the simple representation
of $\Lambda^{u_v(q)}$ of dimension $s_v(\alpha)$ has endomorphism algebra $K$, and hence
so does $X$.

Thus suppose that $q_v=1$. In this case Lemma~\ref{l:simpcond} shows that $\alpha=\epsilon_v$.
Thus $\alpha$ is a root. Clearly also in this case $\End_{\Lambda^q}(X)=K$, as required.

The proof of Theorem~\ref{t:rigidsimp} is analogous to the argument in \cite[\S 4]{CBad}.
Again we work by induction. Assuming either that $\alpha$ is a positive real root, or
that there is a rigid simple of dimension $\alpha$, we again find a vertex $v$ with $(\alpha,\epsilon_v)>0$,
$v$ loopfree and $s_v(\alpha) < \alpha$.

If $q_v\neq 1$ then Theorem~\ref{t:MC} shows that representations of $\Lambda^q$ of
dimension $\alpha$ correspond to representations of $\Lambda^{u_v(q)}$ of dimension $s_v(\alpha)$.
Moreover the conditions
\begin{itemize}
\item[(i)]
$\alpha$ is a positive real root for $Q$,
\item[(ii)]
$q^\alpha=1$, and
\item[(iii)]
there is no decomposition $\alpha=\beta+\gamma+\dots$ as a sum of two or more positive
roots with $q^\beta=q^\gamma=\dots=1$
\end{itemize}
for $q$ and $\alpha$ correspond to the analogous conditions for $u_v(q)$ and $s_v(\alpha)$.
Here we use (\ref{e:usduality}) and the fact that any positive root $\beta$ with
$q^\beta=1$ remains positive under reflection, since the only positive root which
changes sign is $\epsilon_v$ and $q^{\epsilon_v} = q_v \neq 1$.
By induction we get the theorem for $\alpha$.

Thus suppose that $q_v=1$. In this case Lemma~\ref{l:simpcond} shows that there is a simple representation
if and only if $\alpha=\epsilon_v$ and in this case it is clearly rigid.
On the other hand, conditions (i)--(iii) hold if and only if $\alpha=\epsilon_v$, for if (i) and (ii) hold
and $\alpha\neq \epsilon_v$, then the decomposition
\[
\alpha = s_v(\alpha) + \underbrace{\epsilon_v + \dots + \epsilon_v}_{\text{$(\alpha,\epsilon_v)$ terms}}
\]
shows first that
\[
q^{\alpha} = q^{s_v(\alpha)} q^{\epsilon_v} \dots q^{\epsilon_v}
\]
so that $q^{s_v(\alpha)} = 1$, and then that (iii) fails.

\section{The variety of representations}
\label{s:var}
Henceforth $K$ is an algebraically closed field.
In this section we prove Theorem~\ref{t:var}.
Note that we may if we wish assume that $q^\alpha=1$, for otherwise
$\Rep(\Lambda^q,\alpha)$ is empty.
We define $\Rep(L_Q,\alpha)$ to be the open subset of $\Rep(\overline{Q},\alpha)$
consisting of the representations $X$ satisfying (\ref{c:invtwoprime}).
Since it is defined by the nonvanishing of the function
\[
X \mapsto \prod_{a\in \overline{Q}} \det (1_{X_{h(a)}} + X_a X_{a^*}),
\]
we clearly have the following.

\begin{lem}
$\Rep(L_Q,\alpha)$ is a nonempty affine open subset of $\Rep(\overline{Q},\alpha)$.
\end{lem}

There is a morphism of varieties
\[
\nu:\Rep(L_Q,\alpha)\to \GL(\alpha),
\quad
X \mapsto
\biggl(
\prod_{\substack{a\in\overline{Q} \\ h(a)=v}}
(1_{X_v}+X_{a} X_{a^*})^{\epsilon(a)}
\biggr)_{v\in I}.
\]
If $q\in (K^*)^I$, we identify $q$ with the element $(q_v 1_{X_v})_{v\in I}$ of $\GL(\alpha)$, and
then we have $\Rep(\Lambda^q,\alpha) = \nu^{-1}(q)$.
Note that by (\ref{e:detformula}) the image of $\nu$ is actually contained in the subgroup
\[
G = \{ \theta\in\GL(\alpha) \mid \prod_{v\in I} \det(\theta_v) = 1 \}.
\]

\begin{lem}
$\Rep(\Lambda^q,\alpha)$ is an affine variety, and every irreducible component has
dimension at least $g+2p(\alpha)$, where $g = -1+\sum_{v\in I}\alpha_v^2$.
\end{lem}

\begin{proof}
The first statement is clear. The second holds by dimension
theory, since
\[
\dim \Rep(L_Q,\alpha)-\dim G = \sum_{a\in Q} 2\alpha_{h(a)}\alpha_{t(a)} - (-1+\sum_{v\in I} \alpha_v^2)
= g+2p(\alpha)
\]
\end{proof}

Given $X\in \Rep(L_Q,\alpha)$,
we can identify the tangent space $T_X \Rep(L_Q,\alpha)$ with $\Rep(\overline{Q},\alpha)$ and
$T_{\nu(X)} \GL(\alpha)$ with
\[
\End(\alpha) = \bigoplus_{v\in I} \Mat(\alpha_v,K).
\]
The map that $\nu$ induces on tangent spaces is then
\[
d\nu_X : \Rep(\overline{Q},\alpha)\to\End(\alpha),
\quad
d\nu_X(Y) = \sum_{\substack{a\in \overline{Q} \\ h(a)=v}} \ell_a \Delta_a r_a
\]
where
\[
\ell_a = \prod_{\substack{b\in \overline{Q} \\ h(b)=h(a) \\ b<a}}
(1_{X_{h(a)}}+X_b X_{b^*})^{\epsilon(b)},
\quad
r_a = \prod_{\substack{b\in \overline{Q} \\ h(b)=h(a) \\ b>a}}
(1_{X_{h(a)}}+X_b X_{b^*})^{\epsilon(b)},
\]
and
\[
\Delta_a = \begin{cases}
Y_a X_{a^*} + X_a Y_{a^*}
& \text{($\epsilon(a)=1$)} \\
- (1_{X_{h(a)}}+X_a X_{a^*})^{-1} (Y_a X_{a^*} + X_a Y_{a^*}) (1_{X_{h(a)}}+X_a X_{a^*})^{-1}
& \text{($\epsilon(a)=-1$).}
\end{cases}
\]
Now the trace pairing enables one to identify
$\End(\alpha)^*\cong\End(\alpha)$,
and, exchanging the components corresponding to arrows $a$ and $a^*$,
also $\Rep(\overline{Q},\alpha)^* \cong \Rep(\overline{Q},\alpha)$.
Thus the dual of $d\nu_a$ gives a linear map
\[
\phi: \End(\alpha)\to \Rep(\overline{Q},\alpha)
\]
with $\phi(\theta)_a$ equal to
\[
r_a \theta_{h(a)} \ell_a X_a - X_a (1_{X_{t(a)}}+X_{a^*}X_a)^{-1} r_{a^*} \theta_{t(a)} \ell_{a*} (1_{X_{t(a)}}+X_{a^*}X_a)^{-1}
\]
if $\epsilon(a)=1$, and
\[
X_a r_{a^*} \theta_{t(a)} \ell_{a^*} - (1_{X_{h(a)}} + X_a X_{a^*})^{-1} r_a \theta_{h(a)} \ell_a (1_{X_{h(a)}}+X_a X_{a^*})^{-1} X_a
\]
if $\epsilon(a)=-1$.
Now if $\theta$ is in the kernel of $\phi$, then using the identities
\[
X_a (1_{X_{t(a)}} + X_{a^*} X_a) = (1_{X_{h(a)}} + X_a X_{a^*}) X_a
\]
and
\[
\ell_a (1_{X_{h(a)}} + X_a X_{a^*})^{\epsilon(a)} r_a = q_{h(a)} 1_{X_{h(a)}}
\]
one deduces that
\[
q_{h(a)}\theta_{h(a)} \, \ell_a X_a r_{a^*} = \ell_a X_a r_{a^*} \, q_{t(a)}\theta_{t(a)}
\]
Thus by Lemma \ref{l:funnyonto}, $(q_v \theta_v)_{v\in I}$ is an endomorphism of $X$.

If $X$ has trivial endomorphism algebra, then $\phi$ has 1-dimensional kernel.
Thus the image of $d\nu_X$ has codimension 1.
Thus, considering $\nu$ as a map from $\Rep(L_Q,\alpha)$ to $G$,
the induced map on tangent spaces is surjective at $X$. Thus this map defines a
smooth morphism from the set of representations with trivial endomorphism
algebra to $G$. Taking the fibre at $q$, we see that the set of
representations in $\Rep(\Lambda^q,\alpha)$ with trivial endomorphism algebra is smooth.
This completes the proof of Theorem~\ref{t:var}.

\section{Representation type}
\label{s:reptype}
Henceforth we assume that $K$ is algebraically closed of characteristic zero.
The affine quotient scheme $\Rep(\Lambda^q,\alpha)\dbslash \GL(\alpha)$
classifies semisimple representations of $\Lambda^q$ of dimension $\alpha$.
By \cite{LBP} it is stratified into locally closed subsets according to
`representation type'. Here, we say that a semisimple representation $X$
has \emph{representation type} $\tau=(k_1,\beta^{(1)};\dots;k_r,\beta^{(r)})$
provided that it can be decomposed into simple components as
$X = X_1^{\oplus k_1} \oplus \dots \oplus X_r^{\oplus k_r}$
where the $X_i$ are non-isomorphic simples of dimensions $\beta^{(i)}$.
Observe that $\tau$ can only occur if there are indeed simples of dimensions
$\beta^{(i)}$, and that although the $\beta^{(i)}$ need not be distinct,
any real root can occur at most once amongst the $\beta^{(i)}$.
Theorem~\ref{t:var} easily implies the following result. (See for example
the proof of \cite[Theorem 1.3]{CBmm}, except we do not claim irreducibility here.)

\begin{lem}
\label{l:taudim}
If $\tau$ occurs as a representation type for $\Lambda^q$, then the set of
semisimple representations of type $\tau$ is a locally closed subset of
$\Rep(\Lambda^q,\alpha)\dbslash \GL(\alpha)$ of dimension $\sum_{i=1}^r 2p(\beta^{(i)})$.
\end{lem}

Let $\pi:\Rep(\Lambda^q,\alpha)\to \Rep(\Lambda^q,\alpha)\dbslash \GL(\alpha)$
be the quotient map. Observe that the arguments of \cite[\S 6]{CBnorm} can be
applied to the multiplicative preprojective algebra, using
Theorem~\ref{t:homform} instead of \cite{CBex} in the proof of
\cite[Lemma 6.2]{CBnorm}. We obtain the following analogue of \cite[Theorem 6.3]{CBnorm}.

\begin{thm}
If $x$ is an element of $\Rep(\Lambda^q,\alpha)\dbslash \GL(\alpha)$ which has representation type
$(k_1,\beta^{(1)};\dots;k_r,\beta^{(r)})$, then
\[
\dim \pi^{-1}(x) \le g + p(\alpha) - \sum_t p(\beta^{(t)})
\]
where $g = -1+\sum_{v\in I}\alpha_v^2$.
\end{thm}

Combining this with Lemma~\ref{l:taudim}, we have the following.

\begin{cor}
\label{c:invim}
The inverse image in $\Rep(\Lambda^q,\alpha)$ of the stratum of
representation type $(k_1,\beta^{(1)};\dots;k_r,\beta^{(r)})$
has dimension at most $g+p(\alpha)+\sum_{t=1}^r p(\beta^{(t)})$.
\end{cor}

Theorem~\ref{t:simpcrit} now follows.
The variety $\Rep(\Lambda^q,\alpha)$ is equidimensional of dimension $g+2p(\alpha)$
since by Theorem~\ref{t:var} each
irreducible component of $\Rep(\Lambda^q,\alpha)$ has dimension at
least $g+2p(\alpha)$, but by Corollary~\ref{c:invim} and the hypotheses of the theorem,
it has has dimension at most $g+2p(\alpha)$.
Moreover, the non-simple representations form a subset of strictly smaller dimension,
so that the set $\mathcal{S}$ of simple representations must be dense.

\section{Deligne-Simpson problem}
\label{s:DSP}
Fix conjugacy classes $C_1,\dots,C_k$ in $\GL_n(K)$, positive integers
$w=(w_1,\dots,w_k)$ and elements $\xi_{ij}\in K^*$
as in the introduction, so satisfying (\ref{e:xieq}).
Let $Q_w$ and $\alpha$ be the corresponding quiver and dimension vector.
We denote by $a_{ij}$ the arrow with tail at $[i,j]$.
Let $<$ be an ordering on $\overline{Q_w}$ with
$a_{11} < a_{21} < \dots < a_{k1}$.
Define $q\in (K^*)^I$ by $q_0 = 1/\prod_{i=1}^k \xi_{i1}$ and $q_{ij} = \xi_{ij}/\xi_{i,j+1}$.
Observe that $\xi^{[\beta]} = 1/q^\beta$ for any $\beta\in\Z^I$.

\begin{lem}
\label{l:impliesinv}
$\Lambda^q$ is isomorphic to $K\overline{Q_w}/J$ where $J$ is the
ideal generated by the relations
\[
(e_0 + a_{11} a_{11}^*)\dots(e_0 + a_{k1}a_{k1}^*) = q_0 e_0
\]
and
\[
q_{ij} (e_{[i,j]} + a_{ij}^* a_{ij}) = \begin{cases} e_{[i,j]} +
a_{i,j+1}a_{i,j+1}^* & \text{(if $j<w_i-1$)}
\\
e_{[i,j]} & \text{(if $j=w_i-1$).}
\end{cases}
\]
Thus, if $X$ is a representation of $\overline{Q_w}$, then it is a
representation of $\Lambda^q$ if and only if
\[
(1_{X_0} + X_{a_{11}} X_{a_{11}^*})\dots(1_{X_0} + X_{a_{k1}}
X_{a_{k1}^*}) = q_0 1_{X_0}
\]
and
\[
q_{ij} (1_{X_{[i,j]}} + X_{a_{ij}^*}X_{a_{ij}}) = \begin{cases}
1_{X_{[i,j]}} + X_{a_{i,j+1}}X_{a_{i,j+1}^*} & \text{(if
$j<w_i-1$)}
\\
1_{X_{[i,j]}} & \text{(if $j=w_i-1$).}
\end{cases}
\]
\end{lem}

\begin{proof}
It suffices to show that the relations automatically imply (\ref{c:inv}),
for then they are equivalent to (\ref{c:rel}).
For $j=w_i-1$, the relations imply that
$1 + a_{ij}^* a_{ij}$ is invertible in $K\overline{Q_w}/J$. Then
$1+a_{ij} a_{ij}^*$ is invertible (with inverse $1 - a_{ij}(1 +
a_{ij}^* a_{ij})^{-1} a_{ij}^*$), and hence so is $1 + a_{i,j-1}^*
a_{i,j-1}$. Repeating in this way, a descending induction on $j$
gives (\ref{c:inv}).
\end{proof}

We say that a representation $X$ of $\Lambda^q$ is
\emph{strict} if the linear maps $X_{a_{ij}}$
are all injective and the maps $X_{a_{ij}^*}$ are all surjective

\begin{lem}
\label{l:solcorresprep}
There is a representation of $\Lambda^q$ of dimension $\alpha$ if and only if
there are matrices $A_i$ in the conjugacy class closures $\overline{C_i}$ with
$A_1\dots A_k=1$. There is a strict representation $X$ if and only if
there are matrices $A_i\in C_i$ with $A_1\dots A_k=1$.
\end{lem}

\begin{proof}
Given $A_1,\dots,A_k\in \GL_n(K)$, as in \cite[Theorem 2.1]{CBipb} the $A_i$ are in the
closures $\overline{C_i}$ if and only if
there are vector spaces $V_{ij}$ of dimension $\alpha_{ij}$ and linear maps $\phi_{ij}$, $\psi_{ij}$,
\[
K^n \rightoverleft^{\phi_{i1}}_{\psi_{i1}} V_{i1} \rightoverleft^{\phi_{i2}}_{\psi_{i2}} V_{i2}
\rightoverleft^{\phi_{i3}}_{\psi_{i3}} \dots
\rightoverleft^{\phi_{i,w_i-1}}_{\psi_{i,w_i-1}} V_{i,w_i-1}
\]
satisfying
\[
\begin{split}
A_i - \psi_{i1} \phi_{i1}
&= \xi_{i1} \, 1
\\
\phi_{ij} \psi_{ij} - \psi_{i,j+1}\phi_{i,j+1}
&= (\xi_{i,j+1}-\xi_{ij})\, 1_{V_{ij}}
\qquad (1\le j < w_i-1)
\\
\phi_{i,w_i-1} \psi_{i,w_i-1}
&= (\xi_{i,w_i} - \xi_{i,w_i-1})\, 1_{V_{i,w_i-1}}
\end{split}
\]
Moreover, $A_i\in C_i$ if and only if the $\phi_{ij}$ are surjective
and the $\psi_{ij}$ are injective.
Now the equation $A_1\dots A_k=1$ is equivalent to
\[
(\xi_{11} 1+\psi_{11}\phi_{11})\dots(\xi_{k1} 1+\psi_{k1}\phi_{k1}) = 1.
\]
Let $X$ be the representation of $\overline{Q_w}$ defined by
$X_0 = K^n$, $X_{[i,j]} = V_{ij}$, $X_{a_{ij}} = \psi_{ij}/\xi_{ij}$
and $X_{a_{ij}^*} = \phi_{ij}$.
When written in terms of $X$, the equations above are equivalent
to the ones in Lemma~\ref{l:impliesinv}. The lemma follows.
\end{proof}

\begin{lem}
\label{l:irrsolcorresprep}
There is a simple representation of $\Lambda^q$ of dimension $\alpha$ if and only if
there is an irreducible solution to the equation $A_1\dots A_k=1$ with $A_i\in C_i$.
\end{lem}

\begin{proof}
Suppose first that there is an irreducible solution to the equation.
Let $X$ be a strict representation of $\Lambda^q$ corresponding to the solution.
Suppose $Y$ is a subrepresentation of $X$.
Since $A_i = \xi_{i1} (1+X_{a_{i1}}X_{a_{i1}^*})$ we have $A_i(Y_0) \subseteq Y_0$.
Thus irreducibility implies that $Y_0=0$ or $Y_0 = X_0$.
Now if $Y_0 = 0$ then $Y=0$ since the maps $X_{a_{ij}}$ are
all injective, and if $Y_0=X_0$ then $Y=X$ since the maps $X_{a_{ij}^*}$ are
all surjective.
Thus $X$ is simple.

Conversely, suppose that $X$ is a simple representation of $\Lambda^q$ of dimension
vector $\alpha$. We show that $X$ is strict.
If $X_{a_{i\ell}}$ is not injective, let $x\in X_{[i,\ell]}$ be a nonzero element in its kernel.
We define elements $x_j\in X_{[i,j]}$ for $j\ge \ell$ by setting $x_\ell=x$
and $x_{j+1} = X_{a_{i,j+1}^*}(x_j)$ for $j\ge\ell$. An induction using the relation
\[
\xi_{ij} (1_{X_{[i,j]}} + X_{a_{ij}^*}X_{a_{ij}}) = \xi_{i,j+1} (1_{X_{[i,j]}} + X_{a_{i,j+1}}X_{a_{i,j+1}^*})
\]
for $j<w_i-1$ shows that $X_{a_{i,j+1}}(x_{j+1})$ is a multiple of $x_j$ for $j\ge \ell$.
Thus the $x_j$ span a subrepresentation $Y$ of $X$. By simplicity we must have $Y=X$,
but this is impossible since $Y_0=0$ and $\dim X_0 = \alpha_0 = n \neq 0$.
Thus $X_{a_{i\ell}}$ is injective.
A dual argument shows that $X_{a_{i\ell}^*}$ is surjective.
Thus $X$ is strict.
Now let $A_i = \xi_{i1}(1+X_{a_{i1}} X_{a_{i1}^*})$ be the corresponding solution
to $A_1\dots A_k=1$. We show that it is irreducible. Suppose that $Y_0$ is
a subspace of $X_0$ which is invariant under the $A_i$.
Define $Y_{[i,j]}$ inductively by $Y_{[i,1]} = X_{a_{i1}^*}(Y_0)$ and
$Y_{[i,j]} = X_{a_{ij}^*}(Y_{[i,j-1]})$ for $j>1$.
We have
\[
X_{a_{i1}}(Y_{[i,1]}) = X_{a_{i1}} X_{a_{i1}^*} (Y_0) = \frac{1}{\xi_{i1}}(A_i-\xi_{i1}1)(Y_0) \subseteq Y_0,
\]
and then by induction $X_{a_{ij}}(Y_{[i,j]}) \subseteq Y_{[i,j-1]}$ for $j>1$.
Thus $Y$ is a subrepresentation of $X$. By simplicity $Y=0$ or $Y=X$. Thus $Y_0=0$ or $Y_0=X_0$.
Thus the solution is irreducible.
\end{proof}

The proof of Theorem~\ref{t:DSP} is now straightforward.
By \cite[Theorem 1.3]{CBipb} there is a solution to $A_1\dots A_k=1$ with $A_i\in\overline{C_i}$.
By Lemma~\ref{l:solcorresprep} there is a representation of $\Lambda^q$ of dimension $\alpha$.
Now by Theorem~\ref{t:simpcrit} there is a simple representation of $\Lambda^q$
of dimension $\alpha$, and hence by Lemma~\ref{l:irrsolcorresprep} there is
an irreducible solution to $A_1\dots A_k=1$ with $A_i\in C_i$.

\end{document}